\documentclass[draft,12pt]{amsart}

\usepackage{amssymb}
\usepackage{amscd}

\numberwithin{equation}{section}

\theoremstyle{plain}
\newtheorem*{theor}{Theorem 0.1}
\newtheorem{theorem}[equation]{Theorem}

\newtheorem{lemma}[equation]{Lemma}
\newtheorem{corollary}[equation]{Corollary}

\theoremstyle{definition}
\newtheorem{definition}[equation]{Definition}
\newtheorem{example}[equation]{Example}
\newtheorem{examples}[equation]{Examples}
\newtheorem{coexample}[equation]{Counterexample}
\newtheorem{mainexample}[equation]{Main example}
\newtheorem{remark}[equation]{Remark}

\theoremstyle{plain}
\newtheorem*{acknow}{Acknowledgment}

\newcommand{\Sph}{\mathbb S}
\newcommand{\R}{\mathbb R}
\newcommand{\N}{\mathbb N}
\newcommand{\Id}{\operatorname{Id}}
\newcommand{\const}{\operatorname{const}}
\newcommand{\Hd}{\operatorname{Hd}}
\newcommand{\FD}{\operatorname{FD}}
\newcommand{\hb}{\operatorname{hb}}
\newcommand{\pr}{\operatorname{pr}}
\newcommand{\HB}{\mathcal{HB}}
\newcommand{\HS}{\mathcal{HS}}

\begin{document}

{\sloppy
\begin{abstract}{
We discuss two different in general natural approaches to the ideal closure and ideal boundary of Busemann nonpositively curved metric space. It is shown that the identity map of the space admits surjective continuation from its coarse ideal closure to the weak one. We consider some situations when these closures coincide, and when they are essentially different. In particular, the singular Minkowski space is studied as flat Busemann space, and some types of its ideal points are described.
}
\end{abstract}

}

\title{Ideal closures of Busemann space and singular Minkowski space}

\author{P.D.Andreev}
\address{Pomor State University, Arkhangelsk, Russia}
\email{andreev@math.pomorsu.ru}
\date{}

%\thanks{2000 \emph{Mathematics Subject Classification.} Primary , .}
\thanks{Supported by RFBR Grant 04-01-00315-a.}

\maketitle

\vspace{-5mm}

%\noindent

\markboth{{\upshape P.D.Andreev}}{{\upshape Ideal boundaries of Busemann spaces}}

%\keywords    test  \endkeywords
%\subjclass  51  \endsubjclass
\vspace{3mm}

\section*{Introduction}

The boundary at infinity of a metric space $X$, i.e. a set of its points at infinity plays important role in solving a number of problems in metric geometry. Some notions of ideal boundary of metric space are defined in several situations. For example, hyperbolic boundary of Gromov hyperbolic spaces is defined in \cite{Hypgr}, the ideal boundary of simply connected Riemannian manifolds of nonpositive curvature in \cite{BGS}. 
This book contains two descriptions of boundary at infinity of simply connected nonpositively curved manifold $X$. On the first hand, $X$, being metric space, can be embedded into the space $C(X)$ of continuous functions on $X$. The first closure of $X$ is the closure of its image in $C(X)$, or more precisely in its factorspace $C^\star(X) = C(X)/\{\operatorname{constants}\}$. Such defined boundary and closure for general metric space is equivalent to the metric boundary and metric closure introduced in \cite{WW}. We use the term "coarse closure" for this construction here. 

On the other hand, there is well-defined relation "to be asymptotic" on the set of all rays of simply connected nonpositively curved Riemannian manifold $X$. This relation is really equivalence and the set of equivalence classes is the ideal boundary in its second definition. It is shown in \cite{BGS} that the two boundaries are the same in the sense that identity map of $X$ may be continued to homeomorphism of two its ideal closures.

This identity of two approaches to closure of the space remains true in the case of $CAT(0)$-spaces, i.e. simply connected nonpositively curved Aleksandrov spaces. General theory of nonpositively curved Aleksandrov spaces is deeply developed in modern geometry, (see \cite{ABN}, \cite{BH}, \cite{Ba} etc.) and is similar to theory of Riemannian manifolds of nonpositive curvature in a number of results.

The subject of presented paper is boundaries of nonpositively curved spaces in the sense of Busemann and in particular of singular Minkowski spaces. Busemann spaces (\cite{Bu1}, \cite{Bo}) are more general then $CAT(0)$-spaces and nonpositivity of curvature is defined here in more weak sense. One of consequences of this weakness is a possibility of difference between two mentioned approaches to ideal closure.

The simplest examples of such difference arise from singular Min\-kow\-ski spaces. Minkowski spaces appear in Finsler geometry as flat Finsler manifolds. Finsler notion means regularity of their metric. Definition of Busemann spaces given in \cite{Bo} and in \cite{H} leads to consideration of singular Minkowski spaces, in which unit sphere is strictly convex symmetric compact $C^0$-hypersurface but is not necessarily smooth and admits existence of parabolic points.

We call the boundary $\partial_c X$ of Busemann space $X$ arising as boundary of its image when $X$ is embedded into $C^\star(X)$ coarse ideal boundary of $X$, and a set $\partial_w X$ of equivalence classes of asymptotic rays in $X$ its weak ideal boundary. Coarse ideal closure $\bar X_c = X \cup \partial_c X$ has surjective projection $\Pr \colon \bar X_c \to \bar X_w$ onto weak closure, which is the continuation of the identity map of $X$ onto itself. 

More precisely, the theorem holds:

\begin{theorem}\label{proj}
Let $X$ be locally compact Busemann space and $\bar X_c = X \cup \partial_c X$ and $\bar X_w = X \cup \partial_w X$ be correspondingly its coarse and weak closures. Then there is continuous surjective map $\Pr \colon \bar X_c \to \bar X_w$, such that $\Pr|_X = \Id_X$.
\end{theorem}

The simplest example when $\Pr$ is not a homeomorphism arises from the singular Minkowski space or affine space equipped with strictly convex but nonsmooth norm in its directing space. The projection $\Pr$ of coarse ideal boundary of singular Minkowski space occurs related with the inverse Gauss map of its unit sphere $\mathcal S$.

{\sloppy

Let the directing space $V^n$ of singular Minkowski space $A^n$ be equipped with Euclidean structure and $\Sph^{n-1}\subset V^n$ be the unit Euclidean sphere. The Gauss image of the point $y_0 \in \mathcal S$ is the set $\mu(y_0)$ of all vectors $\vec \nu \in \Sph^{n-1}$ such that there is support plane of $\mathcal S^{n-1}$ in the point $y_0$ with external normal $\vec \nu$. Since $\mathcal S$ is strictly convex, the inverse relation $\mu^{-1}$ is a map from $\Sph^{n-1}$ to $\mathcal S$.

}

Given geodesic ray $c\colon [0, +\infty) \to X$, the Busemann function generated by $c$ is the function $\beta_c (y) = \lim\limits_{t \to +\infty} (|y c(t)| - t)$. Every Busemann 
function represents some coarse ideal point $\phi = \beta_c$.
We have following description of two types of Busemann functions on $A^n$. 

\begin{theorem}\label{twoBusemann}
Let $A^n$ be singular Minkowski space. Then coarse ideal point $\phi$ is represented by Busemann function in following two situations:
\begin{enumerate}
\item \label{item1} the direction of $\Pr(\phi)$ is regular;
\item \label{item2} $\phi$ is represented as the limiting horofunction of flag-directed sequence of level 1 with directing flag $\mathcal F(\phi) = (x_0, \bar \alpha_1(\phi))$ where the ray $\bar \alpha_1(\phi)$ has singular direction.
\end{enumerate}
\end{theorem}
See sections \ref{flagdir} and \ref{horofunc} for definitions.

\begin{acknow} My great thanks to V.N.Berestovski\v{\i} for stating the problem and collaboration.
\end{acknow}
\section{Preliminaries}

Let $(X,d)$ be a metric space. The distance between two points $x,y\in X$ will be denoted as $d(x,y) := |xy|$. From now on we will assume that $X$ is locally compact.

\emph{A geodesic} in $X$ is locally isometric immersion $c\colon (a,b) \to X$ of real interval  $(a,b)$, i.e. such immersion that every point $t \in (a,b)$ has a neighbourhood, which is embedded to $X$ isometrically by $c$. Geodesic $c$ is called \emph{minimizing} if  $c$ is isometrical embedding. $c$ is called \emph{complete geodesic} if  $(a,b) = \R$. The interval $[a,b]$ may occur a closed segment, in that case we say that geodesic $c$ \emph{connects points} $y = c(a)$ and $z = c(b)$. The image $c([a,b])$ is called \emph{geodesic segment} in $X$. 
\emph{A ray} in the space $X$ is a minimizing $c\colon [0, +\infty) \to X$, i.e. isometric embedding of half a line $\R_+$ to $X$. 

Metric space  $X$ is \emph{geodesic} if every its point $x$ has a neighbourhood $O(x)$ where any two points $y,z \in O(x)$ can be connected by geodesic segment in $X$ (not necessarily unique). We say that geodesic metric space $X$ is \emph{geodesically complete} if every its geodesic segment is contained in some complete geodesic (not necessarily unique). The Hopf-Rinow  theorem for metric spaces states that locally compact  geodesically complete space $X$ is complete and proper (finitely compact) metric space. 

\begin{definition}(cf. \cite[p. 1904]{H},)
We say that a path $\gamma \colon [0,1] \to X$ is a geodesic segment parameterized proportionally to arc length, if there exists a geodesic $c\colon [a,b] \to X$ such that $\gamma(t) = c(a + t(b-a))$ for all $t \in [0,1]$.
\end{definition}

\begin{definition}\label{Busemann}
Geodesically complete metric space $X$ is called \emph{Busemann space} if for any two geodesic segments  $\gamma_1,\gamma_2: [0,1] \to X$, parameterized proportionally to arc length and such that $\gamma_1(0) = \gamma_2(0)$  inequality
\[|\gamma_1(t)\gamma_2(t)| \le t|\gamma_1(1)\gamma_2(1)|\]
holds for all $t \in [0,1]$.
\end{definition}

If $X$ is a Busemann space, it evidently is contractible and any two its points are connected with unique geodesic segment. Moreover, the definition \ref{Busemann} imply that the metric of Busemann space is convex: for any two paths $\gamma_1,\gamma_2\colon [0,1] \to X$, parameterized proportionally to arc lengths, the function $d(\gamma_1(t),\gamma_2(t))\colon [0,1] \to \R_+$ is convex (\cite{Bu2}).

\begin{examples}
Every $CAT(0)$-space is a Busemann space. Some simply connected Finsler manifolds with reversible metric and nonpositive flag curvature is a Busemann spaces. In particular, every Minkowski space, i.e. finite dimensional affine space equipped with a Finsler metric invariant under translations is a Busemann space. Moreover, every finite dimensional affine space, equipped with strictly convex (not necessarily smooth) norm in the tangent space, which is invariant under translations, is a Busemann space. We will call such a space \emph{singular Minkowski space} (in particular, usual Minkowski space is a singular Minkowski space without any singularity). It was shown in \cite{Bo} that if $X$ admits the action by isometries of cocompact group $G$, then eiser $X$ is Gromov hyperbolic space or it contains Minkowski plane.
\end{examples}

\section{The weak ideal boundary.}

First we recall the definition of boundary of Busemann space given in \cite{H}.

\begin{definition}\label{weak}
Let $X$ be a Busemann space pointed in $x_0 \in X$. \emph{The boundary} of $X$ is a set
\[\partial_{x_{0}} X = \{c\colon \R_+ \to X|\ c\ \mbox{is a ray with}\ c(0)=x_0\}\]
endowed with the compact-open topology.

It was shown in \cite{H} that the definition above is independent on the choice of marked point $x_0$ hence we have well-defined boundary $\partial_w X = \partial_{x_{0}} X$. We will call it \emph{weak boundary} of the space $X$. 
\end{definition}

We will define the topology on the closure of $X$ with this weak boundary.

\begin{definition}\label{defweak}
{\sloppy

We say that sequence $\{\gamma_k\}_{k=1}^{+\infty}$ of length parameterizations $\gamma_k\colon [0, |x_0x_k|] \to \bar X_w$ of segments or rays $[x_0x_k] \subset \bar X_w$ \emph{converges uniformly on compacts} to the length parameterization $\gamma \colon [0, |x_0 y|] \to \bar X_w$ of segment or ray $[x_0y]$ if for any closed finite segment $[a,b] \subset \R_+$ which is containing in infinitely many intervals $[0, |x_0x_k|]$ sequence of restrictions $\gamma_k|_{[a,b]}$ converges uniformly to restriction $\gamma|_{[a,b]}$. Here we denote $[x_0z]$ --- a segment if $z \in X$ or a ray from $x_0$ to $z \in \partial_w X$ if $z$ is a boundary point: $z\in \partial_w X$. In the last case we will mean $|x_0z| = +\infty$ and consider such a ray as segment connecting points $x_0$ and $z \in \partial_w X$.

}
\end{definition}

\begin{definition}
The union $\bar X_w = X \cup \partial_w X$ will be called \emph{weak closure} of $X$. The topology on $\bar X_w$ is \emph{a topology of uniform convergence on compacts of segments}: a sequence $\{y_k\}_{k=1}^{+\infty} \subset \bar X_w$ converges to a point $y \in \bar X_w$ if the sequence of segments $[x_0y_k]$ equipped with length parameterizations $\gamma_k \colon [0, |x_0x_k|] \to \bar X_w$ converges uniformly on compacts to the segment $[x_0y]$ with its length parameterization $\gamma \colon [0, |x_0y|] \to \bar X_w$. 
\end{definition}

There is another description of the weak boundary as a set of equivalence classes of asymptotic rays. 
Two rays $c,d \colon [0,+\infty) \to X$ are \emph{asymptotic} if their Hausdorff distance is finite:
\[\Hd(c,d) < +\infty.\]
The Hausdorff distance between subsets $A,B \subset X$ is 
\[\Hd(A,B) = \inf\{r|\ B \subset \mathcal N_r(A)\ \mbox{and}\ A \subset \mathcal N_r(B)\}.\]
Here $\mathcal N_r(Y)$ is $r$-neighbourhood of a set $Y\subset X$:
\[\mathcal N_r(Y) = \left\{z \in X| \inf\limits_{y \in Y}|yz| < r\right\}.\]

The relation "to be asymptotic" is equivalence on the set of rays in $X$. Weak boundary $\partial_w X$ of $X$ is a set of equivalence classes of asymptotic rays.
This definition is evidently identic to \ref{weak} (cf. \cite{H}). 

The topology of the boundary $\partial_w X$ induced from $\bar X_w$ admits the basis of neighbourhoods $\mathcal U(t, \varepsilon)$ for a point $\xi = c(+\infty) \in \partial_w X$, endpoint of the ray $c$ beginning in $x_0$, where
\begin{equation}
\mathcal U(t, \varepsilon) = \{\eta = d(+\infty) \in \partial_w X|\, d(0) = x_0,\, |c(t), d(t)| < \varepsilon \}.
\end{equation}

\section{The coarse ideal boundary.}

Let $C(X,\mathbb{R})$ be a topological vector space of continuous $\mathbb{R}$-valued functions on $X$ endowed with topology of uniform convergence on compact sets, and $C^*(X,\mathbb{R})$ be its factorspace by the subspace of constants. The topology of $C(X,\mathbb{R})$ arises from the sequence of supremum norms, where $n$-th item is the supremum-norm, defined by the ball of radius $n$ centered in fixed point $x_0\in X$. In particular, when $X$ is proper, this topology coincides with topology of uniform convergence on bounded sets.

Fix a marked point $x_0\in X$. Let $d_y\in C(X,\mathbb{R})$ be the distance function defined by the point $y\in X$ by $d_y(x) := |xy| - |x_0y|$. Every ball $\mathcal{B}(y,r)$ centered in $y$ is its sublevel set:
\[\mathcal{B}(y,r) = \{x \in X|\; d_y(x) \le r - |x_0y|\}.\]
For our space $X$ we have a map $i\colon X \to C(X,\mathbb{R})$, defined as $i(y) = d_y \in C(X,\mathbb{R})$. $i$ is embedding of $X$ into $C(X,\mathbb{R})$. The superposition $p\circ i$, where $p\colon C(X,\mathbb{R}) \to C^*(X,\mathbb{R})$ is a factorization map, defines embedding of $X$ into $C^*(X,\mathbb{R})$, independent on the choice of marked point $x_0$. From now on we will identify the space $X$ with its image $X_* = (p \circ i)(X) \subset C^*(X,\mathbb{R})$.

\begin{remark}\label{Ascoli}
The family of all distance functions $D \subset C(X, \R)$ is obviously equicontinuous when is restricted to arbitrary bounded set. Hence if $X$ is locally compact, then by Arcela-Ascoly theorem this family is precompact on bounded sets, and pointwise convergence in $C(X, \R)$ of functions from $D$ is sufficient for their uniform convergence on compacts.
\end{remark}

\begin{definition}
\emph{Coarse ideal closure} $\bar X_c$ of space $X$ is by definition the closure of its image  $X_*\subset C^*(X,\mathbb{R})$. The set $\partial_c X = \bar X_c \setminus X_*$ is called as \emph{coarse ideal boundary} of $X$, functions forming it as \emph{horofunctions} and we say that every horofunction \emph{presents a coarse points at infinity}. Two horofunctions $\phi,\psi \in C(X, \mathbb{R})$ present the same coarse point at infinity iff they differs by the constant: $\phi - \psi = \const$. The class of horofunction $\phi$ will be considered as a coarse point at infinity, but for simplicity we identify the horofunction $\phi$ and its class and write $\phi \in \partial_c X$.
\end{definition}

In general every horofunction $\phi \in \partial_c X$ is \emph{generated} by a some sequence $\Phi = \{x_k\}_{k=1}^{+\infty}$ as a limiting function of the sequence of corresponding distance functions $\{d_{x_{k}}\}_{k=1}^{+\infty}$.

\begin{definition}
Let $\phi \in C(x, \R)$ be a horofunction, presenting a coarse point at infinity $\xi = [\phi]_s \in \partial_c X$. Sublevels of $\phi$ are called \emph{horoballs centered in $\xi$} and their boundaries, i.e. levels of $\phi$ \emph{horospheres}. We will use notation
\[\HB(\xi, x_0) = \{y | \phi (y) \le \phi(x_0)\}\]
for horoballs with boundary point $x_0$ and
\[\HS(\xi, x_0) = \{y | \phi (y) = \phi(x_0)\}\]
for corresponding horospheres. This notations are independent on the choice of horofunction $\phi$ representing the point $\xi$.
Also, we will consider \emph{open horoballs}
\[\hb(\xi, x_0) = \HB(\xi, x_0) \setminus \HS(\xi, x_0).\]
\end{definition}

\begin{example}
Let $c:[0, +\infty)\to X$ be a ray. The function 
\[\beta_c(x) = \lim\limits_{t \to \infty} (|yc(t)|-t)\]
is called \emph{Busemann function} corresponding to $c$. Evidently, every Busemann function is horofunction. The inverse is also right in simply connected complete nonpositively curved Riemannian manifolds and $CAT(0)$-spaces.
\end{example}

The following counterexample shows that reverse statement may fault in general Busemann space, in particular in singular Minkowski space.

\begin{coexample}\label{coex}
Let $X=A^2$ be two-dimensional singular Min\-kow\-ski space (singular Minkowski plane) with norm 
\begin{equation}\label{singular}
||(y^1,y^2)|| = \sqrt{(y^1)^2 + 2(y^2)^2} + |y^2|
\end{equation}
in the tangent space with coordinates $(y^1, y^2)$. The unit ball centered in the origin in this norm is the intersection of euclidean disks 
$(x^1)^2 + (x^2\pm 1)^2 \le 2$ in coordinates $(x^1, x^2)$ on the plane. Their boundary cycles intersects in two points $(\pm 1, 0)$ orthogonally. Functions  $\epsilon_1 (x^1 +\epsilon_2 x^2) + C$, where $C = \const$ and $\epsilon_1, \epsilon_2 = \pm 1$, are horofunctions as limits of Busemann functions of type 
\[\beta_{\lambda,\mu}(x^1,x^2) = \epsilon_1 (\lambda x^1 + (\mu +\epsilon_2)x^2) + C.\]
Namely, $\beta_{\lambda,\mu} \to \epsilon_1 (x^1 +\epsilon_2 x^2) + C$ when $\lambda \to 1$ and $\mu \to 0$ with condition $\lambda^2 + (\mu +\epsilon_2)^2 =2$.
However they are not Busemann functions theirself, because there is no appropriate ray in $A^2$ to define them as corresponding Busemann functions.
\end{coexample}

The counterexample \ref{coex} shows that the origin point $x_0$ can be connected by the ray not with every coarse point at infinity.

\begin{remark}
In fact, our notion of the coarse ideal closure and coarse ideal boundary is exactly identic to the metric closure and metric boundary described in \cite{WW}, but we act with constructive approach of \cite{BGS} to this metric boundary. It seems to be more convenient for description of asymptotic geometry of Busemann spaces.
\end{remark}

\section{Projection of coarse boundary to the weak one}

Here we will prove the Theorem \ref{proj}. In other words we construct a surjective projection map from the coarse ideal boundary $\partial_c X$ to the weak ideal boundary $\partial_w X$ which maps every Busemann functions to classes of corresponding rays. However, this projection will not be injective in general, even in the set of Busemann functions. First, we will prove the lemma.

\begin{lemma}\label{fouritems}
Let the space $X$ be pointed in $x_0\in X$ and $\phi \in \partial_c X$ be the horofunction with $\phi(x_0)=0$. Then
\begin{enumerate}
\item The function $\phi$ is bounded on any ball in $X$;
\item There exists unique ray $c:[0, +\infty) \to X$ beginning at $c(0)= x_0$ with following property. For every $t>0$ the point $c(t)$ is the unique point of the ball $\mathcal B(x_0, t)$ where $\phi|_{\mathcal{B}(x_0, t)}$ attains its minimum and 
\[\min\phi|_{\mathcal{B}(x_0, t)} =  - t;\]
\item The weak ideal point $\xi \in \partial_w X$ represented by the ray $c$ does not depend on the choice of marked point $x_0$.
\item If $\phi(y) = \lim\limits_{t \to +\infty} (|yc(t)| - t)$ is a Busemann function defined by the ray $c: \R_+ \to X$, then its restrictions to balls $\mathcal B(c(0), t)$ attain their minima in points $c(t)$
\end{enumerate} 
\end{lemma}

\begin{proof}

The first claim of the lemma is the consequence of finitely compactness of the space $X$ and continuity of horofunctions. 

Every distance function $d_y(x) = |xy|- |x_0y|$ with $|x_0y| > r$ when restricted to the ball $\mathcal B(x_0, r)$ attains its unique minimum in the point where the segment $[x_0y]$ exits the ball. The minimal value is
\[\min d_y|_{\mathcal B(x_0, r)} = -r.\]
Hence the minimal value of limiting function $\phi = \lim_{k \to +\infty}$ for a sequence of distance functions $d_{x_{k}}$  is $-r$ and is attained in the boundary sphere.

Let $\phi$ be defined by
\[\phi (y) = \lim\limits_{k \to +\infty} d_{x_{k}}(y),\] 
with the sequence of distance functions $d_{x_{k}}$. Consider the sequence of segments $[x_0x_k]$ and sequence of positive numbers $\varepsilon_i \to 0$. Pick arbitrary $t>0$. For every $i$ there exists $K(i)$ such that for $k>K(i)$ for every $y \in \mathcal B(x_0,t)$
\[|d_{x_{k}}(y) - \phi(y)| < \epsilon_i.\]
In particular, if $z_k\in \mathcal S(x_0, t)$ is exiting point of segment $[x_0x_k]$ with $d_{x_{k}}(z_k) = -t$, then 
\[\phi(z_k) > -t + \varepsilon_i.\]
Hence we have equality for limiting point $z$ of the sequence $z_k$ (it does exist and is unique):
\[\phi(z)=-t.\]
Since segments $[x_0z_k]$ equipped with length parameterizations, converge to the segment $[x_0z]$ and for $p_k \in [x_0z_k]$ with $|x_0p_k| = s \le t$ we have $d_{x_{k}}(p_k) = -s$, then for limiting point $p \in [x_0z]$ with $|x_0p| = s$ equality $\phi(p) = - s$ holds. Since $p$ is the unique point of the ball $\mathcal B(x_0, s)$ with $\phi(p) = -s = \min \phi|_{\mathcal B(x_0,s)}$, then we have well-defined one-parameter family of segments $[x_0z(t)]$ of lengths $|x_0z(t)| = t$ with $[x_0z(s)] \subset [x_0z(t)]$ when $s \le t$ and $\phi(z(t)) = -t$ for all $t > 0$. Their union is the required ray.

Let $x_0'\in X$ be another marked point and $c, c'\colon \R_+ \to X$ be rays beginning in $x_0, x_0'$ correspondingly, such that $c(t)$ and $c'(t)$ are minimal points of the horofunction $\phi$ in balls $\mathcal B(x_0,t)$ and $\mathcal B(x_0', y)$. We prove that rays $c$ and $c'$ are asymptotic. 

Pick $t > 0$ and a sequence $\{x_k\}_{k=1}^{+\infty}$ such that $\phi = \lim\limits_{k \to +\infty} d_{x_{k}}$. Denote $p_k(t)$ a point of the segment $[x_0x_k]$ with $|x_0p_k(t)| = t$ and $p_k'(t)$ a point of segment $[x_0'x_k]$ with $|x_0'p_k'(t)| = t\cdot |x_0'x_k|/|x_0x_k|$. Since the sequence $\{x_k\}_{k=1}^{+\infty}$ is running to infinity, the fraction $|x_0'x_k|/|x_0x_k|$ tends to unity and both sequences $\{p_k(t)\}_{k=1}^{+\infty}$ and $\{p_k'(t)\}_{k=1}^{+\infty}$ converges to points $c(t)$ and $c'(t)$ correspondingly. We have
\[|p_k(t)p_k'(t)| \le \frac{|x_0x_k| - t}{|x_0x_k|}\cdot |x_0x_0'| < |x_0x_0'|\]
and
\[|c(t)c'(t)| \le |x_0x_0'|\]
for all $t$. Hence, Hausdorff distance between rays $c$ and $c'$ is finite:
\[\Hd(c,c') < +\infty.\]

Finally, assume that $\phi(y) = \lim\limits_{t \to +\infty}(|yc(t)| - t)$ is a Busemann function defined by the ray $c\colon \R_+ \to X$. Obviously $\phi(c(t)) = - t$.
\end{proof}

Now we apply the lemma \ref{fouritems} to prove the theorem \ref{proj}

\begin{theor}
Let $X$ be locally compact Busemann space and $\bar X_c = X \cup \partial_c X$ and $\bar X_w = X \cup \partial_w X$ be correspondingly its coarse and weak closures. Then there is continuous surjective map $\Pr \colon \bar X_c \to \bar X_w$ such that $\Pr|_X = \Id_X$.
\end{theor}

\begin{proof}
We put $\Pr(x) = x$ for $x \in X$ and $\Pr(\phi) = c(+\infty)$ for $\phi \in \partial_c X$ with $\phi(x_0)=0$, where $c:[0, +\infty) \to X$ is the ray with $c(0) = x_0$ such that for all $t > 0$
\[\phi(c(t)) = -t = \min\phi|_{\mathcal B(x_0, t)}.\]
This ray is defined by the lemma \ref{fouritems}. The map $\Pr$ is surjective, since for any $\xi \in \partial_w X$
\[\xi = \Pr(\beta_\xi),\]
where $\beta_\xi$ is the Busemann function corresponding to the ray $[x_0\xi]$.

It remains to prove the continuity of the map $\Pr$. Let $\{x_k\}_{k=1}^{+\infty} \subset X$ be the sequence converging to $\phi \in \partial_c X$ in the coarse sense. Then for all $y \in X$
\[d_{x_{k}}(y) \to \phi(y)\]
and for all $t > 0$ restrictions of segments $[x_0x_k]$ to the subsegments of length $t$, converge to the segment $[x_0c(t)]$ of the ray $c$ defined in the lemma \ref{fouritems}. Hence 
\[\Pr(x_k) \stackrel{k\to+\infty}{\longrightarrow} c(+\infty) = \Pr(\phi)\]
in the sense of topology of weak closure $\bar X_w$.
If $\{\phi_k\}_{k=1}^{+\infty}\subset \partial_c X$ is the sequence of horofunctions converging in the coarse sense to the horofunction $\phi \in \partial_c X$, then each its item is the limit of some sequence $\{x_{km}\}_{m=1}^{+\infty}$:
\[\phi_k = \lim\limits_{m \to +\infty} x_{km},\]
and we may choose mixed sequence $\{x_{km(k)}\}_{k=1}^{+\infty}$ converging to $\phi$. Here $m(k)$ is sufficiently large integer defined for all $k$. Since $\Pr(x_{mk}) \to \Pr(\phi_k)$ when $m \to +\infty$, and $\Pr(x_{km(k)}) \to \phi$ when $k \to +\infty$, then $\Pr(\phi_k) \to \Pr(\phi)$ in the weak sense.
\end{proof}

\begin{definition}
Let $\xi \in \partial_c X$ be a coarse point at infinity, represented by the horofunction $\phi$. \emph{The weak projection} $\Pr(\xi)$ of $\xi$ is a weak point at infinity represented by the ray $c\colon \R_+ \to X$ such that for every $t > 0$ the point $c(t)$ is the minimal point of the function $\phi$ restricted onto the ball $\mathcal B(c(0), t)$.  
\end{definition}

\begin{definition}
The Busemann space $X$ is called \emph{regular} if the projection $\Pr$ constructed in the theorem \ref{proj} is homeomorphism of the coarse ideal boundary $\partial_c X$ to the weak one $\partial_w X$. Otherwise it is \emph{singular}.
\end{definition}

\begin{example}
Consider singular Minkowskian plane $A^2$ with coordinates $(x^1,x^2; y^1,y^2)$ in its tangent bundle $TA^2$, where $x(x^1,x^2)$ are coordinates of point $x$ and $\vec y(y^1, y^2)$ --- coordinates of vector $\vec y \in T_xA^2$ which are invariant under translations of $A^2$. The norm in $TA^2$ is defined by \eqref{singular}.

We describe preimage $\Pr^{-1}(\xi)$ where $\xi \in \partial_w A^2$ is an infinite point defined by positively directed $x^1$-axis. The obvious class of horofunctions defining a point in $\Pr^{-1}(\xi)$ is the class of Busemann function $\beta_0$ of the ray $[x_0\xi]$, where $x_0(0,0)$ is the origin. The direct computation gives
\[\beta_0 (x^1,x^2) = |x^2| - x^1.\]
This implies that axes codirected with $[x_0\xi]$ define different classes of horofunctions: the ray $[x_0'\xi]$ with $x_0'(0, x^2_0)$ generates Busemann function 
\[\beta_0' (x^1, x^2) = |x^2-x^2_0| - x^1 \ne \beta_0 (x^1, x^2) + \const.\]
Hence we have one-parameter infinite family of type $\beta_0'$ depending on the second coordinate $x^2_0$.

Two more coarse ideal points which are generated by no ray as Busemann function are mentioned in counterexample \ref{coex}. This horofunctions are
\[\phi_{\pm} (x^1, x^2) = \pm x^2 - x^1.\]
The first of them is generated as the limit of distance functions
\[\phi_+ = \lim\limits_{t \to +\infty} d_{(f(t), -t)}\]
where $f(t)$ is arbitrary function such that 
\[\lim\limits_{t \to +\infty} \frac{f(t)}{t} = +\infty.\]
The second is the limit
\[\phi_- = \lim\limits_{t \to +\infty} d_{(f(t), t)}.\]
\end{example}

\section{Inverse Gauss map in singular Minkowski space.}

Assume that singular Minkowski space $A^n$ is equipped with additional Euclidean structure: Euclidean scalar product $\langle \vec v, \vec w\rangle$ for two arbitrary vectors $\vec v, \vec w \in V^n$. This Euclidean structure has not any relation with Minkowski norm in $V^n$. The Euclidean norm of vector $\vec v \in V^n$ is denoted as 
\[|\,\vec v\,| = \langle \vec v, \vec v \rangle^{1/2}.\]
The unit Euclidean sphere in $A^n$ centered in the origin $x_0$ and unit Euclidean sphere in $V^n$ are denoted as $\Sph^{n-1}$. We will use standard notation $\vec v \perp \vec w$ for perpendicularity of vectors $\vec v$ and $\vec w$ in the Euclidean sense.

The Gauss image $\mu$ of surface $\mathcal S$ is defined for any its point as a set of external normal vectors for support hyperplanes to $\mathcal S$ in this point meant as a set of points of $\Sph^{n-1}$. In general the Gauss correlation is not the map from $\mathcal S$ to $\Sph^{n-1}$, but the inverse correlation is one. The inverse Gauss map
\begin{equation}\label{inverseGauss}
\mu^{-1} \colon \Sph^{n-1} \rightarrow \mathcal S
\end{equation}
sends Euclidean unit vector $\vec \nu$ to the unique point $\mu^{-1}(\vec \nu) \in \mathcal S$ with support hyperplane $\alpha$ such that $\vec \nu$ is its external normal.

Next we study several functions related with the inverse Gauss map.
Fix a hyperplane $A_0^{n-1}$ and its Euclidean normal $\vec \nu_0$. $A_0^{n-1}$ is parallel to support hyperplane to $\mathcal S$ in $\mu^{-1}(\vec \nu_0)$. Define the function $\theta: (\mathcal S \times \mathcal S)_+ \to \R$ by
\[(\vec v - \theta(\vec v, \vec w)\vec w) \perp \vec \nu_0.\]
Here
\[(\mathcal S \times \mathcal S)_+ = \{(\vec v, \vec w)\, | \quad \langle \vec \nu_0, \vec v\rangle > 0, \, \langle \vec \nu_0, \vec w\rangle > 0\}.\]
$\theta(\vec v, \vec w)$ is a norm of projection of $\vec v$ to the $\vec w$-direction in $A_0^{n-1}$-direction. Set
\begin{equation}\label{Theta}
\Theta(\vec v, \vec w) = 
	\left\{
		\begin{array}{cl}
			0& \mbox{if $\vec v = \vec w$}\\
		\frac{1 - \theta(\vec v, \vec w)}{\|\vec v - \theta(\vec v, \vec w)\vec w\|}& \mbox{otherwise}
	\end{array}\right.
\end{equation}

\begin{lemma}\label{Thetalemma}
Let the vector $\vec v_0 = \mu^{-1}(\vec \nu)$ goes in the regular direction. Then the function $\Theta$ is continuous in the point $(\vec v_0, \vec v_0) \in (\mathcal S \times \mathcal S)_+$.
\end{lemma}

\begin{proof}
Let $\vec \nu_i = \vec \nu_0 + \delta_i \vec\nu \in \Sph^{n-1}, \ i = \overline{1,2},\ \vec v_1 = \vec v_0 + \delta_i \vec v = \mu^{-1}(\vec \nu_i) \in \mathcal S$\ with $(\mu^{-1}(\vec \nu_1) \ne \mu^{-1}(\vec \nu_2)$. The elementary computation for Euclidean triangles gives
\begin{equation}\label{cosines}
\frac{(1 - \theta(\vec v_1, \vec v_2)) |\vec v_2|}{|\vec v_1 - \theta(\vec v_1, \vec v_2)\vec v_2|} = \frac{\cos \measuredangle(\vec \nu, \vec v_2 - \vec v_1)}{\cos \measuredangle (\vec \nu, \vec v_2)}.
\end{equation}

All angles here are assumed as Euclidean angles.
The denominator of the right-hand fraction is bounded from zero when the pair $(\vec v_1, \vec v_2)$ belongs to sufficiently small neighbourhood of $(\vec v_0, \vec v_0) \in (\mathcal S \times \mathcal S)_+$, and its numerator is equal to
\[\cos \measuredangle(\vec \nu, \vec v_2 - \vec v_1) = \cos \measuredangle (\vec \nu, \frac{\vec v_2 - \vec v_1}{\|\vec v_2 - \vec v_1\|}).\]
The last item tends to the cosine of the angle between the vector $\vec \nu = \mu(\vec v_0)$ and the limiting plane $A_0^{n-1}$ containing all limiting vectors for $\frac{\vec v_2 - \vec v}{\|\vec v_2 - \vec v_1\|}$, i.e. tangential plane to $\mathcal S$ in the end of $\vec v_0$. The limiting Euclidean angle $\measuredangle(\vec \nu, A_0^{n-1}) = \frac{\pi}{2}$ and
\[\frac{\cos \measuredangle(\vec v_1, \vec v_2 - \vec v_1)}{\cos \measuredangle (\vec \nu, \vec v_2)} \rightarrow 0\]
when $(\vec v_1, \vec v_2) \to (\vec v_0, \vec v_0)$ staying different.
The left-hand fraction of \eqref{cosines} differs from the fraction
\[\Theta(\vec v_1, \vec v_2) = \frac{1 - \theta(\vec v_1, \vec v_2)}{\|\vec v_1 - \theta(\vec v_1, \vec v_2)\vec v_2\|}
\]
only by scaling on bounded value: the correlation $|\vec v|/\|\vec v\|$ is bounded on the set of nonzero vectors since the unit sphere $\mathcal S$ is convex and compact.
Hence
\[\Theta(\vec v_1, \vec v_2) \rightarrow 0\]
when $(\vec \nu, \vec v) \to (\mu(\vec v_0), \vec v_0)$ in any way.
\end{proof}

Consider the set
\begin{equation}
(\Sph^{n-1} \times \mathcal S)_+ = \{(\vec \nu,\vec v)|\, |\vec \nu| = ||\vec v|| = 1, \, \langle \vec \nu, \vec v\rangle > 0\}
\end{equation}
and the continuous function
\begin{equation}
\lambda \colon (\Sph^{n-1} \times \mathcal S)_+ \rightarrow [1, +\infty)
\end{equation}
defined by
\[\vec \nu \perp (\lambda(\vec \nu, \vec v)\cdot \vec v - \mu^{-1}(\vec \nu)).
\]

The value $\lambda(\vec \nu, \vec v)$ is the distance from $x_0$ to intersection point of corresponding support hyperplane to $\mathcal S$ in the point $\mu^{-1}(\vec \nu)$ with the ray eliminating from $x_0$ in $\vec v$-direction.

Set
\begin{equation}
\Lambda(\vec \nu, \vec v) = \left\{
	\begin{array}{cl}
		0& \mbox{if $\mu^{-1}(\vec\nu) = \vec v$}\\
		\frac{\lambda(\vec \nu, \vec v) - 1}{\|\mu^{-1}(\vec\nu) - \vec v\|}& \mbox{otherwise}
	\end{array}\right.
\end{equation}
and for unit vector $\vec v_0 \in \mathcal S$

\begin{equation}\label{L}
L_{\vec v_{0}}(\vec \nu) = \Lambda(\vec \nu, \vec v_0).
\end{equation}

The proofs of both items of following lemma are similar to that of the lemma \ref{Thetalemma} and we omit them.

\begin{lemma}\label{twolambda}
\begin{enumerate}
\item \label{one}
Let the unit vector $\vec v_0 \in \mathcal S$ has regular direction. Then the function $\Lambda$ is continuous in the point $(\mu(\vec v_0), \vec v_0) \in (\Sph^{n-1} \times \mathcal S)_+$;
\item \label{two}
For any vector $\vec v_0 \in \mathcal S$ the function $L_{\vec v_{0}}$ is continuous on the set
\[\Sph_+^{n-1}(\vec v_0) = \{\vec \nu \,|\quad |\vec \nu| = 1, \, \langle \vec \nu, \vec v_0\rangle > 0\}.\]
\end{enumerate}
\end{lemma}

\begin{coexample}
Function $\Theta$ and $\Lambda$ are not continuous in points of type $(\vec \nu, \mu^{-1}(\vec \nu))$, where the vector $\mu^{-1}(\vec \nu))$ goes in singular direction. This can be easily seen in consideration of the counterexample \ref{coex}: for singular vector $\vec v_0 = (1, 0)$ and any external normal $\vec \nu_0$ to its support direction one may find a sequence $\vec v_k \to \vec v_0$ such that $\Lambda(\vec \nu_0, \vec v_k)$ does not converge to $0$. For example if $\vec \nu_0 = (\frac{\sqrt{2}}{2}, \frac{\sqrt{2}}{2})$, such sequence is $\{(\sqrt{2}\cos(\frac{\pi}{4} + \frac1k), 1 - \sqrt{2}\sin(\frac{\pi}{4} + \frac1k))\}_{k=1}^{+\infty}$, if $\vec \nu_0 = (\frac{\sqrt{2}}{2}, -\frac{\sqrt{2}}{2})$ the sequence is $\{(\sqrt{2}\cos(\frac{\pi}{4} + \frac1k), -1 + \sqrt{2}\sin(\frac{\pi}{4} + \frac1k))\}_{k=1}^{+\infty}$, and for all another normals for support lines in $\vec v_0$ one may take both mentioned sequences. 
\end{coexample}

\section{Flag-directed sequences}
\label{flagdir}

It is convenient for us to use following notion of flag in affine space $A^n$.
\begin{definition}
$k$-flag in the space $A^n$ is a $(k+1)$-tiple 
\[\mathcal F^k = (x_0, \bar \alpha^1, \dots, \bar \alpha ^k),\] 
where $\bar \alpha^1$ is a ray emanating from $x_0$ and contained in the line $\alpha^1$, and for all $i \in \overline{2,k}$\ $\bar \alpha^i$ is $i$-dimensional half-plane bounded by $(i-1)$-dimensional plane $\alpha^{i-1}$ and contained in $i$-dimensional plane $\alpha^i$.
\end{definition}

\begin{definition}\label{firstlevel}
A sequence 
\begin{equation} \label{phi}
\Phi = \{x_k\}_{k=1}^{+\infty} \subset A^n
\end{equation}
is called \emph{almost flag-directed sequence of level 1} if the sequence of segments $[x_0x_k]$ converges to a ray $\bar \alpha^1 = [x_0\xi_0]$, where 
\[\xi_0 = \lim\limits_{k \to +\infty} x_k \in \partial_w A^n\]
is a weak limit point at infinity of $\Phi$. The 1-flag $\mathcal F^1 = (x_0, \bar \alpha^1)$ is called its \emph{directing flag}.
\end{definition}

The asymptotic behavior of the sequence \eqref{phi} satisfying definition \ref{firstlevel} can be described by two statements.

\begin{enumerate}
\item For any hyperplane $A^{n-1}_1$ transversal to a line $\alpha^1$ 
\begin{equation}\label{asym}
\lim\limits_{k \to +\infty} d(x_k,\beta) = +\infty.
\end{equation} 
\item If $x_{1,k}$ is a projection of $x_k$ to $A^{n-1}_1$ in $\alpha^1$-direction, then
\begin{equation}\label{rel}
\lim\limits_{k \to +\infty} \frac{|x_0x_{1,k}|}{|x_{1,k}x_k|} = 0.
\end{equation}
\end{enumerate}

\begin{definition}\label{secondlevel}
Let sequence \eqref{phi} be almost flag-directed sequence of level $1$ with directing flag $\mathcal F^1 = (x_0, \bar \alpha^1)$ and $A^{n-1}_1$ be a hyperplane transversal to the line $\alpha^1$.  The sequence $\Phi$ is called \emph{flag-directed sequence of level 1} if the sequence 
\begin{equation}\label{phi1}
\Phi_1 = \{x_{1,k}\}_{k=1}^{+\infty} \subset A^{n-1}_1
\end{equation}
of projections of points of $\Phi$ to $A^{n-1}_1$ is converging to a point 
\[x_{1,0} = \lim\limits_{k \to +\infty} x_{1,k}.\]
The line $\alpha_{1,0}$ passing throw $x_{1,0}$ in $\alpha_1$-direction is a \emph{asymptotic line} of the sequence $\Phi$. 

$\Phi$ is called  \emph{almost flag-directed sequence of level 2}, if the sequence  \eqref{phi1} is almost flag-directed sequence of level 1 with directing ray $\bar \alpha_{1,1} \subset A^{n-1}_1$. \emph{The directing flag} of almost flag-directed sequence $\{x_k\}_{k=1}^{+\infty}$ of level 2 is 2-flag $\mathcal F^2 = (x_0, \bar \alpha_1, \bar\alpha_2)$, where 2-dimensional half-plain $\bar\alpha_2$ is bounded by the line $\alpha_1$ and contains the ray $\bar\alpha_{1,1}$.
\end{definition}

The definition \ref{secondlevel} is independent on the choice of the hyperplane $A^{n-1}_1$ transversal to $\alpha^1$. Next we will continue this construction recurrently to define flag-directed sequences of arbitrary level $p \le n$. Assume that we have already defined the notion of flag-directed sequences of levels up to $p-1$ and almost flag-directed surfaces of levels up to $p$ and their directing flags. The directing flag $\mathcal F^{p-1} = (x_0, \bar\alpha_1, \dots \bar\alpha_{p-1})$ of flag-directed sequence $\Phi=\{x_k\}_{k=1}^{+\infty}$ includes $(p-1)$ half-planes of dimensions from 1 to $(p-1)$, and corresponding flag 
\begin{equation}\label{p-flag}
\mathcal F^p = (x_0, \bar\alpha_1, \dots \bar\alpha_p)
\end{equation}
for almost flag-directed sequence of level $p$ includes half-planes of dimensions from 1 to $p$.

\begin{definition}
Let $\Phi$ as in \eqref{phi} be almost flag-directed sequence of level $p-1$ with directing flag \eqref{p-flag}. Let $A^{n-p}_p$ be $(n-p)$-plane passing throw $x_0$ transversal to $\alpha_p$ and $x_{p,k}$ be the projection of the point $x_k$ to $A^{n-p}_p$ in $\alpha_p$ direction. The sequence $\Phi$ is called \emph{flag-directed of level $p$} with \emph{directing flag \eqref{p-flag}} if the sequence 
\[\Phi_p=\{x_{p,k}\}_{k=1}^{+\infty}\subset A^{n-p}\]
 is converging to a point 
\[x_{p,0} = \lim\limits_{k \to +\infty} x_{p,k}.\]
The $p$-plane $\alpha_{p,0}$ passing throw $x_{p,0}$ in $\alpha_p$-direction is a \emph{asymptotic plane} of the sequence $\Phi$. 

The sequence $\Phi$ is \emph{almost flag-directed of level $(p+1)$} if $\Phi_p$ is almost flag-directed of level 1 with directing ray $\bar\alpha_{1,p} \subset A^{n-p}_p$. Its \emph{directing flag} is 
 \[\mathcal F^{p+1} = (x_0, \bar\alpha_1, \dots, \bar\alpha_p, \bar\alpha_{p+1})\] 
where $\bar\alpha_{p+1}$ is $(p+1)$-dimensional half-plane bounded by $p$-plane $\alpha_p$ and containing the ray $\bar\alpha_{1,p}$. Every almost flag directed sequence of level $n$ is defined \emph{flag-directed of level $n$}. \emph{The asymptotic plane} of such sequence is the whole space $A^n$. 
\end{definition}

This definition again is independent on the choice of the $(n-p)$-plain $A^{n-p}_p$ transversal to $\alpha_p$. The theorem \ref{extract} below claims that every sequence of points in $A^n$ admits picking out the subsequence which is converging or flag-directed. Hence from precompactness property of the set of distance functions, every coarse ideal point is the limiting point of some flag-directed sequence.

\begin{mainexample}
Let $(x^1, \dots, x^n)$ be the affine coordinate system in $A^n$. Every sequence of points with coordinates 
\begin{equation}\label{main}
x_k = (f_1(k), \dots, f_p(k), g_1(k), \dots, g_{n-p}(k)),
\end{equation}
where $\lim\limits_{k \to +\infty}f_i(k) = +\infty$, functions $g_j(k)$ converges to a finite limit when $k \to +\infty$ and for all $i \in \{1, \dots, p-1\}$
\[\lim\limits_{k \to + \infty}\frac{f_{i+1}(k)}{f_i(k)} = 0.\]
Then the sequence \eqref{main} is flag-directed sequence of level $p$. Moreover, every flag-directed sequence of level $p$ admits the representation \eqref{main} in some affine coordinate system.
\end{mainexample}

The following theorem and its corollary reveals the purpose of introduction of flag-directed sequences.

\begin{theorem}\label{extract}
Given sequence $\Psi = \{x_k\}_{k=1}^{+\infty}$ of points in $A^n$ one may extract the subsequence which is converging or flag-directed. 
\end{theorem}

\begin{proof}
Assume that $\Psi$ contains no bounded subsequence.

Consider the sequence \[\left\{\frac{\overrightarrow{x_0y_m}}{|x_0y_m|}\right\}_{m=1}^{+\infty} \subset \mathcal S\]
of unit vectors in $y_m$-directions from $x_0$. By compactness of $\mathcal S$ it has converging subsequence corresponding to subsequence 
\[\Psi _1 = \{y_{m_{i_{1}}}\}_{i_1=1}^{+\infty} \subset \{y_m\}_{m=1}^{+\infty}.\] 
Set $\bar \alpha_1$ be the directing ray of
\[\vec v_1 = \lim\limits_{i_1 \to +\infty} \frac{\overrightarrow{x_0y_{m_{i_{1}}}}}{|x_0y_{m_{i_{1}}}|}.\]
The subsequence $\Psi$ is almost flag-directed of level 1. If it is flag-directed, we finish the procedure. If not, the following step is to choose its almost flag-directed subsequence of level 2 and so on. The procedure will be finished in at most $n$ steps. We will describe its second step and $p$-th step.

Consider the sequence of unit bivectors 
\begin{equation}\label{bivectors}
\left\{\frac{\vec v_1 \wedge \vec w_{m_{i_{1}}}}{\|\vec v_1 \wedge \vec w_{m_{i_{1}}}\|}\right\}_{i_{1}=1}^{+\infty}  \subset G_+(2,n),
\end{equation}
where $\vec w_{m_{i_{1}}}$ is the directing vector of 2-dimensional half-plane bounded by the line $\alpha_1 \supset \bar \alpha_1$ and containing the point $y_{m_{i_{1}}}$.  Here $G_+(2,n)$ is Grassmannian manifold of oriented two-dimensional planes in $A^n$ represented as the set of its unit bivectors. The norm of bivector may be considered in the sense of satellite Euclidean structure on $A^n$.

From compactness of $G_+(2,n)$ the sequence \eqref{bivectors} contains converging subsequence corresponding to a subsequence
\[\Psi_2 = \left\{y_{m_{{{i_{1}}_{i_{2}}}}}\right\}_{i_{2}=1}^{+\infty}.\]
The subsequence $\Psi_2$ is almost flag-directed of level 2 and its directing flag is \[\mathcal F_2 = (x_0, \bar \alpha_1, \bar \alpha_2)\]
where $\bar \alpha_2$ is the limiting half-plane directed by the limiting bivector of the sequence \eqref{bivectors}. If the subsequence $\Psi_2$ is flag-directed, we finish the procedure.

Let we have already constructed the almost flag-directed subsequence of level $p-1$
\[\Psi_{p-1} = \left\{y_{m_{j_{p-1}}}\right\}_{j_{p-1}=1}^{+\infty} \subset \Psi.\] 
Its directing flag is 
\[\mathcal F_{p-1} = (x_0, \bar \alpha_1, \dots, \bar \alpha_{p-1}).\]
Denote as $\vec Q$ the directing $(p-1)$-vector of the oriented $(p-1)$-plane $\alpha_{p-1}$.

Consider the sequence of unit $p$-vectors 
\begin{equation}\label{polivectors}
\left\{\frac{\vec Q \wedge \vec w_{m_{j_{p-1}}}}{\|\vec Q \wedge \vec w_{m_{j_{p-1}}}\|}\right\}_{j_{p-1}=1}^{+\infty}  \subset G_+(p,n),
\end{equation}
constructed similarly to \eqref{bivectors}. By the compactness of the Grassmannian manifold $G_+(p,n)$ it has converging subsequence and hence we may choose almost flag-directed subsequence of $\Psi$ of level $p$. Since every almost flag-directed sequence of level $n$ is flag-directed, the procedure will be finished by finite number of steps with extracting of flag-directed subsequence from $\Psi$.
\end{proof}

The proof of the corollary below is based on the theorem \ref{convergence} which is proved independently.

\begin{corollary}
Given horofunction $\phi \in \partial_c A^n$ there exists flag-directed sequence $\Phi = \{x_k\}_{k=1}^{+\infty}$ such that
\[\phi = \lim\limits_{k\to +\infty}d_{x_{k}}.\]
\end{corollary}

\begin{proof}
Every horofunction is by definition a limit of distance functions. Let the sequence $\Psi = \{y_m\}_{m=1}^{+\infty}$ satisfy $\phi = \lim\limits_{m \to +\infty} d_{y_{m}}$. If $\Psi$ is flag-directed, the theorem is proved. Otherwise we will extract the flag-directed subsequence $x_k = y_{m_{k}}$ from $\Psi$. Application of the theorem \ref{convergence} concludes the proof.
\end{proof}

\begin{definition}
\emph{The level} of the horofunction $\phi \in \partial_c A^n$ is the minimal level of the flag-directed sequence $\Phi \in \pr^{-1} (\phi) \subset \FD$.
The directing flag of such minimal level sequence $\Phi$ is called \emph{directing flag} of the horofunction $\phi$.
\end{definition}

Our next purpose is to establish some relations between coarse ideal points of singular Minkowski space $A^n$ and flag-directed sequences.

\section{Defining horofunctions via flag-directed sequences.}
\label{horofunc}

Here we will show that each flag-directed surface generates a horofunction as a limit function of corresponding sequence of distance functions. The space $A^n$ is the singular Minkowski space. First we pay our attention on two simplest cases. Following lemma gives the proof of item \ref{item2} of the theorem \ref{twoBusemann}.

\begin{lemma}\label{firstlevellemma}
Let \eqref{phi} be flag-directed sequence of level 1. Then the sequence
\begin{equation}\label{distphi}
d_\Phi = \{d_{x_{k}}\}\subset C(A^n)
\end{equation}
converges to a horofunction which is Busemann function.
\end{lemma}

\begin{proof}
The flag $\mathcal F^1 = (x_0, \bar \alpha_1)$ corresponding to the sequence $\Phi$ defines the direction of the ray $\bar \alpha_1$ and the weak ideal point $\xi = \bar\alpha_1(+\infty)$. Let $\beta_1\subset A^n$ be arbitrary hyperplane transversal to $\bar\alpha_1$ and
\[x_{1,0} = \lim\limits_{k \to +\infty} x_{1,k}\]
where $x_{1,k}$ is the projection of $x_k$ to $\beta_1$ in $\alpha_1$ direction. Segments $[x_{1,0}x_k]$ converges to the ray $[x_{1,0}\xi]$ codirected with $\bar\alpha_1$. Set $c\colon [0, +\infty) \to A^n$ the arc length parameterization of $[x_{1,0}\xi]$. This defines Busemann function \[\beta_{(c, x_{1,0})}(y) = \lim\limits_{t \to +\infty}(|yc(t)|-t).\]
Evidently we may choose the hyperplane $A^{n-1}_1$ such that $\beta_{(c, x_{1,0})}(x_0)=0$. We claim that $\beta_{(c, x_{1,0})}$ is the proposed horofunction. Really, for arbitrary $y \in A^n$
\[|d_{x_k}(y) - \beta_{(c,x_{1,0})}(y)| = ||x_ky| - |x_kx_0| - \beta_{(c,x_{1,0})}(y)|.
\]
If $t_k \in [0, +\infty)$ is such a number that vector $\overrightarrow{c(t_k)x_k}$ is parallel to $A^{n-1}_k$, then
\[\lim\limits_{k \to +\infty} |c(t_k)x_k| = 0\]
and for any $\varepsilon > 0$ there exists $K_1 \in \N$, such that if $k > K_1$, then
\[|c(t_k)x_k| < \frac{\varepsilon}{3}\]
and hence
\begin{equation}\label{001}
||x_ky| - |c(t_k)y||< \frac{\varepsilon}{3}.
\end{equation}
By definition of Busemann function there exists $K_2 \in \N$ such that for all $k > K_2$
\begin{equation}\label{002}
||c(t_k)y| - t_k - \beta_{(c,x_{0,1})}(y) | < \frac{\varepsilon}{3}
\end{equation}
and $K_3 \in \N$ such that for all $k > K_3$
\begin{equation}\label{003}
|t_k - |x_0c(t_k)|| < \frac{\varepsilon}{3}
\end{equation}
Combining estimations \eqref{001}, \eqref{002} and \eqref{003} one gets for $k > K = \max \{K_1, K_2, K_3\}$
\[|d_{x_{k}}(y) - \beta_{(c,x_{0,1})}(y)| < \varepsilon\]
and sequence of functions $d_{x_{k}}$ converges to Busemann function $\beta_{(c,x_{0,1})}$ pointwise. Using the remark \ref{Ascoli} completes the proof.
\end{proof}

Some more general statement is formulated in following lemma.

\begin{lemma}\label{almost1}
Let \eqref{phi} be almost flag-directed surface of level 1 with directing flag $\mathcal F = (x_0, \bar \alpha_1)$. Then if the sequence $d_{x_{k}}$ converges in the point $y_0$, then it converges in any point $y$ such that $(y_0y) \parallel \alpha_1$ and
\[\lim\limits_{k \to +\infty} (d_{x_{k}}(y) - d_{x_{k}}(y_0)) = -(t - t_0)\]
where $t, t_0$ are parameters of points $y, y_0$ correspondingly in the arc length parameterization $\gamma: \R \to A^n$ of the line $(y_0y)$ codirected with $\bar \alpha$.
\end{lemma}

\begin{proof}
The claim can be easily obtained from equality
\[d_{x_{k}}(y) - d_{x_{k}}(y_0) = |x_ky| - |x_ky_0|\]
and converging of both sequences of rays $\{[y_0x_k)\}_{k=1}^{+\infty}$ and $\{[yx_k)\}_{k=1}^{+\infty}$ to rays codirected with $\bar \alpha$.
\end{proof}

This lemma has a number of applications and the first one is the proof of the item \ref{item1} of the theorem \ref{twoBusemann}.

\begin{lemma}\label{usefullemma}
Let $\Phi$ be almost flag-directed sequence of level 1 with directing flag $\mathcal F^1 = (x_0,\bar \alpha_1)$, where the ray $\bar\alpha_1$ has regular direction. Then the sequence \eqref{distphi} converges to Busemann function of the ray $\bar\alpha_1$.
\end{lemma}
\begin{proof}
Regularity of direction for $\alpha_1$ means that the sphere $\mathcal S$ has a tangential plane as support cone at the end of unit vector $\vec v$ codirected with $\bar\alpha_1$. Take a hyperplane $A_1^{n-1}$ to be parallel to this tangent hyperplane. We will show that for all $y_0 \in A_1^{n-1}$
\begin{equation}\label{distzero}
\lim\limits_{k \to +\infty} d_{x_{k}}(y_0) = 0.
\end{equation}
For this we consider the homothety $h_k : A^n \to A^n$ with coefficient  $\kappa_k = |x_0x_k|^{-1}$ which moves $x_k$ to $x_0$. Set $h_k(x_0) = z_k \in \mathcal S$, $h_k(A_1^{n-1}) = \Pi_k$, $h_k(y_0) = y_k \in \Pi_k$, $\vec v_k = \overrightarrow{x_0z_k}$ and $\vec w_k$ be the unit vector, codirected with $\overrightarrow{x_0w_k}$. Then both sequences of vectors $\{\vec v_k\}_{k=1}$ and $\{\vec w_k\}_{k=1}$ converge to the vector $v_0$ contradirected to $\bar \alpha_1$. By the lemma \ref{Thetalemma} the function $\Theta(\vec v, \vec w)$ defined in \eqref{Theta} is continuous in the point $(\vec v_0, \vec v_0)$ and consequently
\[\lim\limits_{k \to +\infty} \frac{1 - |x_0y_k|}{|y_kz_k|} \to  0.\]
Hence for any $\varepsilon > 0$ there exists $K\in \N$, such that for all $k > K$
\[|1 - |x_0y_k|\, | < \frac{\varepsilon}{|x_0y_0|} |y_kz_k|.\] 

Applying the inverse homothety to $h_k$, one gets
\[|d_{x_{k}} (y_0)| = |\, |x_kx_0| - |x_k y_0|\, | = |1 - |x_0y_k|\, | \cdot |x_kx_0| < \varepsilon.\] and equality \eqref{distzero}. Lemma \ref{almost1} gives that the sequence $\Phi$ converges to some horofunction. Next we show that this horofunction coincides with the Busemann function of the ray $\bar \alpha_1$.

Let $c:\R_+ \to A^n$ be the arc length parameterization of the ray $\bar \alpha_1$ and
\[\beta_1(y) = \lim\limits_{t \to +\infty} (|yc(t)| - t|)\]
be corresponding Busemann function.

First, take the point $y_0 \in A_1^{n-1}$.
The homothety $\tilde h_t$ with coefficient $\frac1t$ which moves the point $c(t)$ to $x_0$ translates the hyperplane $A_1^{n-1}$ to tangent hyperplane $\tilde \Pi$ to $\mathcal S$ in the point $\tilde x_0 = \tilde h_t(x_0)$. Denote $\tilde y_t = \tilde h_t(y_0)$, $\vec v_0 = \overrightarrow{x_0\tilde x_0}$, $\vec w_t$ --- unit vector codirected with the ray $[x_0 \tilde y_t)$ and $\vec \nu_t$ --- the external normal to $\mathcal S$ in the end of the vector $\vec w_t$. With this notation the function $L_{\vec v_0}(\vec \nu_t)$ defined in \eqref{L}, tends to $0$  when $t \to +\infty$ (the item \eqref{two} of the lemma \ref{twolambda}), and for any $\varepsilon > 0$ there exists $T > 0$ such that for all $t > T$
\[\frac{\|x_0\tilde y_t\| - 1}{\|\tilde y_t \tilde x_0\|} < \frac{\varepsilon}{\|x_0y_0\|}.\]
Applying again the inverse homothety to $\tilde h_t$, one gets
\[\|y_0c(t)\| - t < \varepsilon\]
and
\[\beta_1(y_0) = \lim\limits_{t \to +\infty}  \|y_0c(t)\| - t = 0 = \lim\limits_{k \to +\infty} d_{x_{k}}(y_0).\]
From the lemma \ref{almost1} the identity holds for all points $y \in A^n$:
\[\beta_1(y) = d_{x_{k}}(y).\]
\end{proof}

Now we are ready  to prove the main theorem of the paragraph.

\begin{theorem}\label{convergence}
For any flag-directed sequence $\Phi = \{x_k\}_{k=1}^{+\infty}$ in the singular Minkow\-ski space $A^n$, the sequence of corresponding distance functions $\{d_k\}_{k=1}^{+\infty}$ converges to some horofunction.
\end{theorem}

\begin{proof}
The base of induction is proved in the lemma \ref{firstlevellemma}

Assume that the statement is true for any flag-directed surface in $A^n$ of level less than $p \le n$ and consider the flag-directed sequence $\Phi$ of level $p > 1$ with directing flag $\mathcal F = (x_0, \bar \alpha_1, \dots, \bar \alpha_p)$. For any $(n-p)$-plane $A_p^{n-p}$ transversal to the $p$-plane $\alpha_p$, the sequence $\{x_{p,k}\}_{k=1}^p$ converges to a point $x_{p,0} = \lim_{k \to +\infty} x_{p,k}$. Let $A_1^{n-1}$ be the hyperplane containing $A_p^{n-p}$ and transversal to the line $\alpha_1$. Then the sequence $\Phi_1 = \{x_{1,k}\}_{k=1}^{+\infty}$ is a flag-directed sequence of level $(p-1)$ in $A_1^{n-1}$ and hence in $A^n$. By assumption, the sequence of corresponding distance functions $d_{x_{1,k}}$ is converging to a horofunction $\phi_1$ defined in $A_1^{n-1}$ and by the lemma \ref{usefullemma} in $A^n$.

Fix a $(n-2)$-plane $\Pi_0$ containing $\alpha_p$, and let $A_{1,k}^{n-1}$ to be the hyperplane containing $\Pi_0$ and the point $x_k$. Since $p>1$, then only finitely many of hyperplanes $A_{1,k}^{n-1}$ are parallel to $\alpha_1$ and when $A_{1,k}^{n-1}$ is transversal to $\alpha_1$, then, by the induction assumption, the sequence $\{x_{1,k,p}\}_{p=1}^{+\infty}$ of projections to $A_{1,k}^{n-1}$ in $\alpha_1$-direction of points $x_p$ generates the limiting horofunction
\[\phi_k = \lim\limits_{p \to +\infty} d_{x_{1,k,p}}.\]
It is not hardly to verify, using methods of the lemmas \ref{distphi} -- \ref{usefullemma} that horofunctions $\phi_k$ converges pointwise, and hence uniformly on bounded sets, to the limiting horofunction $\phi = \lim\limits_{k \to +\infty} \phi_k$ which coincides with the limit
\[\phi = \lim\limits_{k \to +\infty} d_{x_{k}}.\]
the theorem is proved by the induction.
\end{proof}

The proved theorem allows to define the projection map 
\[\pr \colon \FD \rightarrow \partial_c A^n\]
from the set $\FD$ of all flag-directed sequences in $A^n$ to its coarse ideal boundary. For the sequence $\Phi \in \FD$ its image $\pr(\Phi)$ is its limiting horofunction.

We finish the paper with two simple statements describing conditions for two sequences to generate the same horofunction.

\begin{theorem}
If two flag-directed sequences $\Phi_1, \Phi_2 \in \FD$ have the same flag $\mathcal F$ and the same asymptotic plane, then $\pr(\Phi_1) = \pr(\Phi_2)$.
\end{theorem}

\begin{proof}
The sequence $\Phi$ obtained from $\Phi_1$ and $\Phi_2$ with alternating their items is also flag-directed with the same directing flag and asymptotic plane. The horofunction, generated by $\Phi$ coincides with the horofunction generated by any its subsequence, in particular by both $\Phi_1$ and $\Phi_2$.
\end{proof}

Consequently, $\pr(\Phi)$ depends only on directing flag and asymptotic plane of $\Phi \in \FD$.

The following statement is not much more complicated.

\begin{theorem}\label{rigid}
Let two almost flag-directed sequences of level $p$ 
\[\Phi_i = \{x_{i,k}\}_{k=1}^{+\infty}\ i = \overline{1,2}\] have common flag $\mathcal F = (x_0, \bar \alpha_1, \dots, \bar \alpha_p)$ and for all $k$ 
\[\overrightarrow{x_{1,k}x_{2,k}} \parallel \alpha_p.\] 
Then if the sequence $\Phi_1$ generates the horofunction $\phi$, then $\Phi_2$ generates $\phi$ as well.
\end{theorem}

\begin{proof}
It is easy to see that the unique limiting point of the sequence $\{d_{x_{2,k}}\}$ is $\phi \in \partial_c A^n$. Since this sequence is precompact in $C(A^n, \R)$, it converges, and its limit is namely $\phi$.
\end{proof}

\end{document}